\documentclass[12pt]{amsart}

\usepackage{amsthm, amsmath, amssymb, mathrsfs}
\usepackage{enumerate}

\usepackage[colorlinks]{hyperref}

\makeatletter
\@namedef{subjclassname@2010}{%
	\textup{2010} Mathematics Subject Classification}
\makeatother

\allowdisplaybreaks

\theoremstyle{plain}
\newtheorem{thm}{Theorem}[section]
\newtheorem{cor}[thm]{Corollary}
\newtheorem{lem}[thm]{Lemma}
\newtheorem{prop}[thm]{Proposition}
\theoremstyle{definition}
\newtheorem{defn}[thm]{Definition}
\theoremstyle{remark}
\newtheorem{rem}[thm]{Remark}
\numberwithin{equation}{section}
\newtheorem*{nota}{Notation}

\begin{document}
	
	\title[Ideals in the dual of introverted subspaces of $PM_\Psi(G)$]{Ideals in the dual of introverted subspaces of $\Psi$-pseudomeasures} 
	\author[A. Dabra]{Arvish Dabra$^\star$}
	\address{Arvish Dabra, \newline\indent Department of Mathematics, \newline\indent Indian Institute of Technology Delhi, \newline\indent New Delhi - 110016, \newline\indent India.}
	\email{arvishdabra3@gmail.com}
	\author[R. Lal]{Rattan Lal}
	\address{Rattan Lal, \newline\indent Department of Mathematics, \newline\indent
		Punjab Engineering College (Deemed to be University), \newline\indent
		Chandigarh - 160012, \newline\indent India.}
	\email{rattanlal@pec.edu.in}
	
	\author[N. S. Kumar]{N. Shravan Kumar}
	\address{N. Shravan Kumar, \newline\indent Department of Mathematics, \newline\indent
		Indian Institute of Technology Delhi, \newline\indent
		New Delhi - 110016, \newline\indent India.}
	\email{shravankumar.nageswaran@gmail.com}
	
	\keywords{Generalised Fourier algebras; Maximal regular ideals; Minimal ideals; Orlicz spaces; Topologically introverted subspaces.\\
 $\star$ Corresponding author.\\
 \textit{email address:} arvishdabra3@gmail.com (A. Dabra).}
 
\begin{abstract}
Let $G$ be a locally compact group and $(\Phi,\Psi)$ a complementary pair of Young functions satisfying the $\Delta_2$-condition. Let $A_\Phi(G)$ be the Orlicz analogue of the Fig\`{a}-Talamanca Herz algebra $A_p(G).$ The dual of the algebra $A_\Phi(G)$ is the space of $\Psi$-pseudomeasures, denoted by $PM_\Psi(G).$ For certain topologically introverted subspaces $\mathcal{A}$ of $PM_\Psi(G)$ and the Banach algebras $W_\Phi(G)$ or $B_\Phi(G),$ denoted by $\mathcal{B},$ we characterise the maximal regular left/right/two-sided ideals of the Banach algebras $\mathcal{A}^{'}$ and $\mathcal{B}^{''}$ considered with the Arens product.  We further characterise the minimal left ideals of $\mathcal{A}^{'}$ and prove the necessary and sufficient conditions for the existence of minimal ideals in the algebras $A_\Phi(G)$ and $\mathcal{B}.$
\end{abstract}

\maketitle

\section{Introduction}
The study of the structure of closed left/right/two-sided ideals in several Banach algebras is an intriguing problem in harmonic analysis. Given a locally compact group $G,$ the Group algebra $L^1(G),$ the Measure algebra $M(G),$ the Fourier algebra $A(G),$ the von Neumann algebra $VN(G)$ and the algebra of convolution operators are some of the significant Banach algebras.

It is well known that for any Banach algebra $B,$ the second conjugate space $B^{''}$ when equipped with the Arens product also becomes a Banach algebra. In the case when $B$ is \textit{commutative}, Filali \cite{filali} characterised the maximal regular left/right/two-sided ideals of $B^{''}.$ Furthermore, he also studied the structure of minimal left/right ideals of $B.$ In 2011, Filali and Monfared \cite{FM} characterised all finite-dimensional left ideals in the \textit{dual of introverted subspaces} of $B^{'}.$

 For any locally compact \textit{abelian} group $G,$ Ghahramani and Lau \cite{GL} studied the maximal regular ideals in the second conjugate space of $L^1(G,w),$ where $w$ represents a \textit{weight} function on $G.$ In 1964, Eymard \cite{Eym} introduced and extensively studied the Fourier algebra $A(G)$ as well as its dual $VN(G).$ The algebra $A(G)$ is isometrically isomorphic to $L^1(\widehat{G}),$ where $\widehat{G}$ is the dual group of an \textit{abelian} group $G.$ Lau \cite{Lau} studied in detail the structure of maximal regular left ideals of $VN(G)^{'},$ which is the second conjugate space of the Fourier algebra $A(G).$ 

For $1 < p < \infty,$ Fig\`{a}-Talamanca \cite{figa} and Herz \cite{herz1971} generalised the algebra $A(G)$ to $A_p(G)$ by replacing $L^2(G)$ with $L^p(G).$ The algebra $A_p(G)$ is a commutative Banach algebra and is referred to as Fig\`{a}-Talamanca Herz algebra. Ghahramani and Lau \cite{GL} studied the structure of maximal regular ideals of $PM_p(G)^{'} = A_p(G)^{''}.$ The Banach algebra $A_p(G)$ has received significant attention from many researchers such as Daws \cite{MD}, Derighetti \cite{DeriS, DeriI, DeriB}, Forrest \cite{FO1,FO2}, Granirer \cite{Gra0,Gra, GraS}, Hosseini and Amini \cite{Amini} and others in the recent decades. Over the past few years, Daws and Spronk \cite{daws2} and Gardella and Thiel \cite{GH1,GH2} have worked on the convolution algebras on $L^p$-spaces and studied the structure of the algebra of $p$-pseudomeasures $PM_p(G).$ The study of ideals in the \textit{dual of certain topologically introverted subspaces} of $PM_p(G) = A_p(G)^{'}$ by Derighetti et al. \cite{DFM} inspires the work of our article.

Since Orlicz spaces are a natural generalisation of the classical $L^p$-spaces, Lal and Kumar \cite{RLSK1} and Aghababa and Akbarbaglu \cite{AA} independently introduced and studied the Orlicz Fig\`{a}-Talamanca Herz algebra, denoted by $A_\Phi(G).$ The algebra $A_\Phi(G)$ is the $L^\Phi$-analogue (a generalisation) of the Fig\`{a}-Talamanca Herz algebra $A_p(G).$ For any locally compact group $G,$ the space $L^\Phi(G)$ refers to the Orlicz space associated with the Young function $\Phi.$ If $\Phi \in \Delta_2,$ then the dual of the algebra $A_\Phi(G)$ is isometrically isomorphic to the space of $\Psi$-pseudomeasures $PM_\Psi(G),$ where $\Psi$ corresponds to the complementary function to $\Phi$ \cite[Theorem 3.5]{RLSK1}.

The main objective of this article is to extend the work of Derighetti et al. \cite{DFM} to the Orlicz analogues and characterise the left/right/two-sided maximal regular ideals and minimal left ideals in the \textit{dual of certain topologically introverted subspaces} of $PM_\Psi(G) = A_\Phi(G)^{'}.$ The article focuses on the following closed subspaces, namely, $C_{\delta,\Psi}(\widehat{G}), PF_\Psi(G), M(\widehat{G}), AP_\Psi(\widehat{G}), WAP_\Psi(\widehat{G}), UCB_\Psi(\widehat{G})$ or $PM_\Psi({G})$ and we denote each of these space by $\mathcal{A}.$

This article is structured in the following manner: In Section \ref{sec2}, we recall the necessary background that is needed for the rest of the article. In Section \ref{sec3}, we characterise the maximal regular left/right/two-sided ideals of the Banach algebras $\mathcal{A}^{'}$ and $\mathcal{B}^{''},$ where $\mathcal{B}$ is either $W_\Phi(G)$ or $B_\Phi(G).$ In Theorem \ref{s3t1}, we prove that if $M$ is any maximal regular left/right/two-sided ideal of $\mathcal{A}^{'},$ then $M$ is either weak$^\ast$ dense in $\mathcal{A}^{'}$ or there exists a unique $x \in G$ such that $M = \{ \Gamma \in \mathcal{A}^{'}: \langle \Gamma, \lambda_\Psi(x) \rangle = 0\}.$ In Section \ref{sec4}, the structure of minimal left ideals of $\mathcal{A}^{'}$ is studied. Theorem \ref{s4t6} establishes that $M$ is a minimal left ideal in $\mathcal{A}^{'}$ if and only if $M = \textit{\textbf{k}} \, \Gamma,$ where $\Gamma$ is either topologically $x$-invariant for some $x \in G$ or a non-zero right annihilator of $\mathcal{A}^{'}.$ In addition, we also demonstrate the existence of minimal ideals in the Banach algebras $A_\Phi(G)$ and $\mathcal{B}.$

\section{Preliminaries}\label{sec2}

Let us begin this section by recalling some fundamentals of the Orlicz space theory.

\subsection{Orlicz space}
 A convex function $ \Phi: \mathbb{R}\rightarrow [0,\infty]$ is termed as a Young function if it is symmetric and satisfies $\Phi(0)= 0$ and $\lim\limits_{t \to \infty} \Phi(t)= + \infty.$ Note that, for a given Young function $\Phi$, the map $\Psi$ defined by $$\Psi(s):= \sup{\{t \, |s|-\Phi(t):t\geq 0\}}, \hspace{1cm} s\in\mathbb{R},$$ is also a Young function. It is called as the complementary function to $\Phi.$ Further, the pair $(\Phi,\Psi)$ (as well as $(\Psi,\Phi)$) is referred to as a complementary pair of Young functions. For $1 < p < \infty,$ the function $\Phi(t) = |t|^p/p$ is an example of a Young function with $\Psi(s) = |s|^q/q$ as its complementary function, where $1/p + 1/q = 1.$

Let $G$ be a locally compact group with a left Haar measure $dt.$ A Young function $\Phi$ is said to satisfy the $\Delta_{2}$-condition, denoted $\Phi\in\Delta_{2},$ if there exists a constant $c > 0$ and $t_{0} > 0$ such that the inequality $\Phi(2t)\leq c \, \Phi(t)$ holds for all $t\geq t_{0}$ whenever $G$ is compact and for non compact $G,$ the same inequality holds with $t_{0}=0.$ If the Young function $\Phi$ and its complementary function $\Psi$ both satisfy the $\Delta_2$-condition, then we call that the pair $(\Phi,\Psi)$ satisfies the $\Delta_2$-condition. It is easy to verify that the family $(\Phi_\beta,\Psi_\beta)$ of complementary pair of Young functions satisfies the $\Delta_2$-condition, where $\Phi_\beta(t) := |t|^\beta \, (1 + |\log|t||)$ for $\beta > 1.$

The Orlicz space associated to a Young function $\Phi,$ denoted $L^{\Phi}(G),$ is defined as $${L}^{\Phi}(G) := \left\{ f: G \rightarrow  \mathbb{C}:f \, \mbox{is measurable and}\int_G\Phi(\alpha |f|)\ dt <\infty \text{ for some}~ \alpha>0  \right\}.$$ The space $L^{\Phi}(G)$ becomes a Banach space with respect to the Luxemburg (or Gauge) norm $N_{\Phi}(\cdot),$ given by $$N_{\Phi}(f):= \inf \left\{a>0:\int_G\Phi\left(\frac{|f|}{a}\right) dt \leq1 \right\}.$$ Let $\Psi$ be the complementary function to $\Phi.$ The Orlicz norm $\|\cdot\|_{\Phi}$ on $L^\Phi(G),$ is defined by $$\|f\|_{\Phi} := \sup \Bigg \{\int_{G}|fg| \, dt :g\in L^\Psi(G) \, \, \text{and} \, \, \int_{G}\Psi(|g|)dt\leq1 \Bigg\}.$$ It is well known that these two norms are equivalent.
If a Young function $\Phi \in \Delta_{2},$ then the space $\mathcal{C}_c(G)$ of all continuous functions on $G$ with compact support is dense in $L^\Phi(G).$ Further, if the pair $(\Phi,\Psi)$ satisfies the $\Delta_2$ condition, then the dual of $(L^{\Phi}(G),N_{\Phi}(\cdot))$ is isometrically isomorphic to $(L^{\Psi}(G),\|\cdot\|_{\Psi}).$ In particular, the space $L^\Phi(G)$ is reflexive.

\noindent We refer \cite{RR} for more information on the Orlicz space.

\subsection{Orlicz Fig\`{a}-Talamanca Herz algebra}
We begin this subsection with the definition of $A_\Phi(G).$ Assume that the pair $(\Phi,\Psi)$ satisfies the $\Delta_2$-condition. For any function $f:G\rightarrow\mathbb{C}$, we define $\check{f}$ as $\check{f}(t):=f(t^{-1})$ for all $t \in G.$ The space $A_\Phi(G)$ is defined as the set of all continuous functions $u$ on $G$ that are of the form $$u=\sum_{n \in \mathbb{N}} f_n \ast \check{g_n},$$ where $f_n\in L^\Phi(G)$ and $g_n\in L^\Psi(G)$ such that $$\sum_{n \in \mathbb{N}} N_\Phi(f_n)\|g_n\|_\Psi<\infty.$$ The space $A_\Phi(G)$ equipped with the norm $$\|u\|_{A_{\Phi}}:= \inf \left\{\sum_{n \in \mathbb{N}} N_{\Phi}(f_{n})\|g_{n}\|_{\Psi}: u=\sum_{n \in \mathbb{N}} f_{n} \ast \check {g}_{n} \right\}$$ is a commutative Banach algebra with pointwise addition and multiplication. This algebra is called Orlicz Fig\`{a}-Talamanca Herz algebra. Furthermore, it is a regular, tauberian, semi-simple Banach algebra with Gelfand spectrum homeomorphic to $G.$

\noindent For more details on the Banach algebra $A_\Phi(G),$ we refer the readers to the series of papers \cite{AA,AAR,AD,RLSK1,RLSK3,RLSK2,RLSK4}.

\subsection{Spaces and algebras related to $A_\Phi(G)$}

We begin by defining the multiplier algebra of $A_\Phi(G)$, denoted $B_\Phi(G).$ The space $B_\Phi(G)$ consists of continuous functions $u$ such that $uv \in A_\Phi(G) \, \, \forall \, v \in A_\Phi(G).$ For any $u \in B_\Phi(G),$ one can define a bounded linear map $T_u: A_\Phi(G) \to A_\Phi(G)$ by $T_u(v) := uv.$ Then $B_\Phi(G)$ is a commutative Banach algebra with the operator norm and pointwise addition and multiplication. It is clear that $A_\Phi(G) \subseteq B_\Phi(G).$

Let $M(G)$ be the set of all bounded complex Radon measures. For $\mu \in M(G)$ and $g \in L^\Psi(G),$ define the map $T_\mu: L^\Psi(G) \to L^\Psi(G)$ by $T_\mu(g) := \mu \ast g.$ If $\mathcal{B}(L^\Psi(G))$ denotes the space of all bounded linear operators on $L^\Psi(G)$ with the operator norm, then it is easy to verify that $T_\mu \in \mathcal{B}(L^\Psi(G)).$ Let $PM_\Psi(G)$ be the closure of the set $\{T_\mu: \mu \in M(G)\}$ in $\mathcal{B}(L^\Psi(G))$ with respect to the ultra-weak topology (weak$^*$-topology). It follows from \cite[Theorem 3.5]{RLSK1} that for a locally compact group $G,$ the dual of $A_\Phi(G)$ is isometrically isomorphic to $PM_\Psi(G).$ Let $PF_\Psi(G)$ be the norm closure of the set $\{T_g: g \in L^1(G)\}$ inside $\mathcal{B}(L^\Psi(G)).$ The dual of $PF_\Psi(G)$ is $W_\Phi(G),$ where $W_\Phi(G)$ is a commutative Banach algebra containing $A_\Phi(G)$ \cite{RLSK4}. For $\mu = \delta_x \, \, (x \in G),$ we denote the operator $T_{\delta_x}$ by $\lambda_\Psi(x)$ and observe that $\lambda_\Psi$ represents the left regular representation of $G$ on $\mathcal{B}(L^\Psi(G)).$ From \cite{AD} and \cite{RLSK3}, let us recall the following subspaces of $PM_\Psi(G)$ that are frequently used in this article.
\begin{align*}
 C_{\delta,\Psi}(\widehat{G}) &:= \overline{span\{\lambda_\Psi(x):x \in G\}}, \\ M(\widehat{G}) &:= \overline{M(G)}, \\ UCB_\Psi(\widehat{G}) &:= \overline{A_\Phi(G)\cdot PM_\Psi(G)}, \\ AP_\Psi(\widehat{G}) &:= \{ T \in PM_\Psi(G): A_\Phi(G) \to PM_\Psi(G), u \mapsto u \cdot T \, \, \text{is compact}\}, \\ WAP_\Psi(\widehat{G}) &:= \{ T \in PM_\Psi(G): A_\Phi(G) \to PM_\Psi(G), u \mapsto u \cdot T \, \, \text{is weakly compact}\}.
\end{align*}
We remark that these subspaces are an Orlicz analogue of the subspaces of $PM_{p}(G)$ which are considered in \cite{Gra}.

\subsection{Banach algebras and Topologically introverted subspaces} Let $B$ be a Banach algebra and $X$ be a closed linear subspace of $B^{'}.$ Then $X$ is said to be \textit{left topologically invariant} if $T \cdot a \in X$ for all $a \in B$ and $T \in X,$ where $(T \cdot a)(b) := T(ab),$ for $b \in B.$ Given such subspace $X$ of $B^{'},$ one can define a continuous linear functional $\Gamma \, \odot \, T$ on $B$ by $$\langle \Gamma \, \odot \, T,a\rangle := \langle \Gamma, T \cdot a \rangle,$$ for $a \in B, T \in X$ and $\Gamma \in X^{'}.$ If $\Gamma \, \odot \, T \in X$ for all $T \in X$ and $\Gamma \in X^{'},$ then $X$ is termed as a \textit{left topologically introverted} subspace of $B^{'}.$ Given a left topologically introverted subspace $X$ of $B^{'},$ it is easy to verify that $X^{'}$ is a Banach algebra with the Arens product given by
$$\langle \tilde{\Gamma} \, \square \, \Gamma, T\rangle := \langle \tilde{\Gamma},\Gamma \odot T \rangle,$$
for $\tilde{\Gamma}, \Gamma \in X^{'}$ and $T \in X.$
Similarly, one can have the notion of \textit{right} topologically introverted subspace. However, in our case, since $B = A_\Phi(G)$ is a commutative Banach algebra, both these notions coincide.

Throughout this article, $G$ denotes a locally compact group with a fixed Haar measure and the pair $(\Phi,\Psi)$ of complementary Young functions satisfies the $\Delta_2$-condition. Recall that we denote the following subspaces
$$C_{\delta,\Psi}(\widehat{G}), PF_\Psi(G), M(\widehat{G}), AP_\Psi(\widehat{G}), WAP_\Psi(\widehat{G}), UCB_\Psi(\widehat{G}) \, \, \text{or} \, \, PM_\Psi({G}),$$
of $PM_\Psi(G)$ by $\mathcal{A},$ and we denote the Banach algebras
$$W_\Phi(G) \, \, \text{or} \, \, B_\Phi(G),$$
by $\mathcal{B}.$

\section{Maximal regular ideals}\label{sec3}

Our main aim in this section is to characterise the maximal regular left/right/two-sided ideals of the Banach algebras $\mathcal{A}^{'}$ and $\mathcal{B}^{''}.$ The ideas and results are inspired by Derighetti et al. \cite{DFM}.

Let us begin with a significant lemma about the space $\mathcal{A}.$

\begin{lem}\label{TI}
	The space $\mathcal{A}$ is a norm closed topologically introverted subspace of $PM_\Psi({G}).$
\end{lem}

\begin{proof}
It is easy to verify that the space $\mathcal{A}$ is a norm closed topologically invariant subspace of $PM_\Psi(G).$ For $\mathcal{A} = UCB_\Psi(\widehat{G}),$ the space is topologically introverted as shown in \cite[Proposition 2.6]{FST1}. Since the rest of the spaces are norm closed invariant subspace of $WAP_\Psi(\widehat{G}),$ the result is a consequence of \cite[Lemma 1.2 (b)]{LL}.
\end{proof}

From the preceding lemma, it follows that the space $\mathcal{A}^{'}$ is a Banach algebra with the Arens product $\square.$ The next result provides an identification of the Banach algebra $A_\Phi(G)$ inside the algebras $\mathcal{A}^{'}$ and $\mathcal{B}^{''}.$ More precisely,

\begin{thm}\label{IWC}
The algebra $A_\Phi(G)$ can be identified with a subalgebra contained in the center of the Banach algebras $\mathcal{A^{'}}$ and $\mathcal{B}^{''}.$
\end{thm}

\begin{proof}
It follows from \cite[Corollary 3.2]{RLSK4} that
$$A_\Phi(G) \subseteq W_\Phi(G) \subseteq B_\Phi(G),$$
i.e., $A_\Phi(G)$ is a subalgebra of the commutative Banach algebra $\mathcal{B}.$ The fact that $A_\Phi(G)$ can be identified with a subalgebra in the center of $\mathcal{B^{''}}$ follows from \cite[Lemma 3.9]{CY}. 

For $\mathcal{A} = PM_\Psi(G),$ the result again follows from \cite[Lemma 3.9]{CY} since $A_\Phi(G)$ is also a commutative Banach algebra. For the case $\mathcal{A} = PF_\Psi(G),$ since $PF_\Psi(G)^{'} = W_\Phi(G)$ is again a commutative Banach algebra, the result is clear. For the rest of the cases, i.e., $\mathcal{A} = C_{\delta,\Psi}(\widehat{G}), M(\widehat{G}), AP_\Psi(\widehat{G}), WAP_\Psi(\widehat{G}), \, \, \text{or} \, \, UCB_\Psi(\widehat{G}),$ the space $\mathcal{A}$ is a norm closed subspace of $PM_\Psi(G).$ By \cite[Theorem 10.1, Pg. 88]{conway}, we have
$$\mathcal{A}^{'} = PM_\Psi(G)^{'} \slash \mathcal{A}^\perp,$$
where $\mathcal{A}^\perp := \{\Gamma \in PM_\Psi(G)^{'}: \Gamma(\mathcal{A}) = \{0\}\}.$ Now, consider the continuous linear map $q \circ J: A_\Phi(G) \to \mathcal{A}^{'},$ where $J$ is the canonical mapping from $A_\Phi(G)$ to $PM_\Psi(G)^{'}$ and $q$ is the quotient map from $PM_\Psi(G)^{'}$ to $PM_\Psi(G)^{'} \slash \mathcal{A}^\perp.$ We denote the image of any element $u \in A_\Phi(G)$ under the map $q \circ J$ with $u$ itself. Since $\mathcal{A}$ contains the set $\{\lambda_\Psi(x):x\in G\}$ and the Gelfand spectrum of $A_\Phi(G)$ is homeomorphic to $G,$ it follows that the map $q \circ J$ is injective. Further, observe that for $u \in A_\Phi(G)$ and $T \in \mathcal{A},$ 
$$|\langle u,T \rangle| = |\langle J(u)+\mathcal{A}^\perp,T\rangle| = |\langle J(u),T\rangle| = |\langle T,u\rangle| \leq \|u\|_{A_\Phi} \, \|T\|,$$
i.e., $\|u\|_{\mathcal{A}^{'}} \leq \|u\|_{A_\Phi}.$ The final conclusion follows from \cite[Lemma 3.9]{CY} and the above-mentioned identification of $A_\Phi(G)$ in $\mathcal{A}^{'}.$
\end{proof}

\begin{rem}\label{s3r1}
    We remark that $A_\Phi(G)$ is weak$^\ast$ dense in $\mathcal{A}^{'}.$ By \cite[Theorem A.3.29 (i)]{dales}, this is clear for $\mathcal{A} = PM_\Psi(G).$ For the rest of the cases, let $\Gamma_1 \in \mathcal{A}^{'}.$ Then $\Gamma_1 = \Gamma + \mathcal{A}^\perp$ for some $\Gamma \in PM_\Psi(G)^{'}.$ Again by  \cite[Theorem A.3.29 (i)]{dales}, there exists a net $\{u_\alpha\}$ in $A_\Phi(G)$ such that $J(u_\alpha) \to \Gamma$ in $\sigma(PM_\Psi^{'},PM_\Psi),$ where $J$ is the canonical mapping from $A_\Phi(G)$ to $PM_\Psi(G)^{'}.$ Now, for every $T \in \mathcal{A},$
    $$\langle u_\alpha,T \rangle = \langle J(u_\alpha)+\mathcal{A}^\perp,T \rangle = \langle J(u_\alpha),T \rangle \to \langle \Gamma,T\rangle = \langle \Gamma + \mathcal{A}^\perp,T \rangle = \langle \Gamma_1,T\rangle,$$
    i.e., $u_\alpha \to \Gamma_1$ in $\sigma(\mathcal{A}^{'},\mathcal{A}).$
\end{rem}

The next corollary is a straightforward consequence of \cite[Proposition 4.5]{DFM} and Theorem \ref{IWC}.

\begin{cor}\label{s3c1}
    Let $M$ be a maximal regular left/right/two-sided ideal in $\mathcal{A}^{'}$ $($or $\mathcal{B}^{''}).$ Then either $M \cap A_\Phi(G) = A_\Phi(G),$ or $M \cap A_\Phi(G) = I_x := \{u \in A_\Phi(G): u(x) = 0\}$ for a unique $x \in G.$
\end{cor}

Recall that \cite[Pg. 101]{Herz} if $B$ is an algebra of functions, then the support of a linear functional $T \in B^{'}$ as a subset of $G$ is characterised by: $x \notin supp(T)$ if and only if there exists a neighbourhood $V$ of $x$ such that $T(v)=0$ for all $v \in B$ with $supp(v) \subseteq V.$ The following proposition gives the structure of the maximal regular ideals of $\mathcal{A}^{'}$ under certain conditions.

\begin{prop}
    Let $M$ be a maximal regular left/right/two-sided ideal in $\mathcal{A}^{'}.$ Assume that there exists a non-zero $T \in \mathcal{A}$ such that $M \subseteq \{\Gamma \in \mathcal{A}^{'}: \langle \Gamma,T\rangle = 0\}.$ Then $M$ equals $M_x,$ where $M_x := \{\Gamma \in \mathcal{A}^{'}: \langle \Gamma,\lambda_\Psi(x)\rangle = 0\}$ for a unique $x \in G.$
\end{prop}

\begin{proof}
    Since $T \in \mathcal{A}$ is non-zero, by Corollary \ref{s3c1}, $$M \cap A_\Phi(G) = I_x = \{u \in A_\Phi(G): u(x) = 0\} \subseteq M$$ for a unique $x \in G.$ By the given hypothesis, $T(u) = 0$ for all $u \in I_x.$ We claim that $supp(T) = \{x\}.$ Let if possible, there exists $y \in supp(T)$ such that $y \neq x.$ Choose $V$ to be the compact neighborhood of $y$ disjoint from $x.$ Since $y \in supp(T),$ there exists $v \in A_\Phi(G)$ with $supp(v) \subseteq V$ and $T(v) \neq 0.$ However, this is absurd as $v \in I_x.$ Thus, our claim follows. Now, by \cite[Theorem 3.6]{RLSK1}, as singletons are sets of spectral synthesis for $A_\Phi(G),$ we have,
    $T = c \, \lambda_\Psi(x)$ for some constant $c.$ This implies that 
    $$M \subseteq \{\Gamma \in \mathcal{A}^{'}: \langle \Gamma,\lambda_\Psi(x)\rangle = 0 \} = M_x.$$
    Since $M_x$ is also an ideal, the result follows from the maximality of $M.$
\end{proof}

Note that in the above proposition, for $\mathcal{A} = PF_\Psi(G),$ the action $\langle \Gamma,\lambda_\Psi(x)\rangle$ is understood as $\Gamma(x).$ This action is well defined as $PF_\Psi(G)^{'} = W_\Phi(G) \subseteq \mathcal{C}(G),$ i.e., $\Gamma$ is a continuous function on $G.$ The need to define this action arises as the set $\{\lambda_\Psi(x):x \in G\}$ may not be contained in $PF_\Psi(G)$ for any arbitrary group $G.$

\begin{nota}
For a Banach algebra $B,$ we denote the Gelfand spectrum of $B$ by $\sigma(B).$
\end{nota}

Here is the promised result of this section that gives the characterisation of the maximal regular ideals of the Banach algebras $\mathcal{A}^{'}$and $\mathcal{B}^{''}.$ This serves as an Orlicz analogue of \cite[Theorem 4.8]{DFM}.

\begin{thm}\label{s3t1}
   $(1)$ Let $M$ be a maximal regular left/right/two-sided ideal of $\mathcal{A^{'}}.$ Then $M$ is either weak$^\ast$ dense in $\mathcal{A^{'}},$ or there exists a unique $x \in G$ such that $$M = M_x =\{\Gamma \in \mathcal{A}^{'}: \langle \Gamma,\lambda_\Psi(x)\rangle = 0\}.$$

   $(2)$ If $M$ is a maximal regular left/right/two-sided ideal of $\mathcal{B}^{''},$ then $M$ is either weak$^\ast$ dense in $\mathcal{B}^{''},$ or $M = M_\varphi :=\{\Gamma \in \mathcal{B}^{''}: \langle \Gamma,\varphi\rangle = 0\}$ for a unique $\varphi \in \sigma(\mathcal{B}).$
\end{thm}

\begin{proof}
$(1)$ Let $M$ be a maximal regular left/right/two-sided ideal of $\mathcal{A}^{'}.$ We denote the weak$^\ast$ closure of $M$ in $\mathcal{A}^{'}$ by $\overline{M}^{w^\ast}.$ Our claim is that $\overline{M}^{w^\ast}$ is also a left/right/two-sided ideal of $\mathcal{A}^{'}.$ Let $\Gamma \in \overline{M}^{w^\ast}$ and $\Gamma_1 \in \mathcal{A}^{'}$ be arbitrary. Then there exists a net $\{\Gamma_\alpha\}$ in $M$ such that $\Gamma = \lim\limits_\alpha \Gamma_\alpha$ and by Remark \ref{s3r1}, there exists a net $\{u_\beta\}$ in $A_\Phi(G)$ such that $\Gamma_1 = \lim\limits_\beta u_\beta.$ By \cite[Theorem 2.6.15 (ii)]{dales}, the product in $\mathcal{A}^{'}$ is weak$^\ast$ continuous on the left and by Theorem \ref{IWC}, since $A_\Phi(G)$ is contained in the center of $\mathcal{A}^{'},$ we have,
$$\Gamma_1 \, \square \, \Gamma = \lim_\beta (u_\beta \, \square \, \Gamma) = \lim_\beta (\Gamma \, \square \, u_\beta) = \lim_\beta \lim_\alpha (\Gamma_\alpha \, \square \, u_\beta) = \lim_\beta \lim_\alpha (u_\beta \, \square \, \Gamma_\alpha) \in \overline{M}^{w^\ast}$$
if $M$ is a left ideal of $\mathcal{A}^{'}.$
If $M$ is a right ideal, then
$$\Gamma \, \square \, \Gamma_1 = \lim_\alpha (\Gamma_\alpha \, \square \, \Gamma_1) \in \overline{M}^{w^\ast}.$$
Hence, the claim is true. Now, $M \subseteq \overline{M}^{w^\ast}$ and since $M$ is a maximal ideal, either $\overline{M}^{w^\ast} = \mathcal{A}^{'}$ or $\overline{M}^{w^\ast} = M.$ In the first case, it follows that $M$ is weak$^\ast$ dense in $\mathcal{A}^{'}.$ Otherwise, by Corollary \ref{s3c1}, $M \cap A_\Phi(G) = I_x = \{u \in A_\Phi(G): \langle \lambda_\Psi(x),u\rangle = u(x) = 0\}$ for a unique $x \in G.$ As $I_x$ is a maximal regular ideal of $A_\Phi(G),$ one can prove by repeating the same arguments as in \cite[Theorem 5.3]{CY} that the weak$^\ast$ closure of $I_x$ in $\mathcal{A}^{'}$ is a maximal ideal of $\mathcal{A}^{'}$ and is given by $\overline{I_x}^{w^\ast} = \{\Gamma \in \mathcal{A}^{'}: \langle \Gamma,\lambda_\Psi(x)\rangle = 0\} = M_x.$ Thus,
$$M_x = \overline{I_x}^{w^\ast} \subseteq \overline{M}^{w^\ast} = M.$$
Now, by maximality of $M_x$ and $M,$ it follows that $M = M_x.$

$(2)$ Since $\mathcal{B}$ is a commutative Banach algebra, the result is a consequence of \cite[Theorem 3.2]{filali}.
\end{proof}

\begin{lem}\label{s3l3}
The Banach algebra $A_\Phi(G)$ is an ideal in $\mathcal{A^{'}}$ if the group $G$ is discrete. The converse holds for $\mathcal{A}$ other than $PF_\Psi(G).$ 
\end{lem}

\begin{proof}
As $G$ is discrete, it follows from \cite[Lemma 4.9]{AD} and \cite[Corollary 4.13 (ii)]{RLSK3} that
$$C_{\delta,\Psi}(\widehat{G}) = PF_\Psi(G) = M(\widehat{G}) = UCB_\Psi(\widehat{G}) \subseteq AP_\Psi(\widehat{G}) \subseteq WAP_\Psi(\widehat{G}) \subseteq PM_\Psi(G).$$ The proof of the lemma follows verbatim to the proof of \cite[Theorem 4.3]{AD} by replacing $PM_\Psi(G)$ with $\mathcal{A}.$ Hence, we omit it.
\end{proof}

\begin{rem}
    In the case when $\mathcal{A} = PF_\Psi(G),$ since its dual $W_\Phi(G)$ is a commutative Banach algebra containing $A_\Phi(G),$ the converse of Lemma \ref{s3l3} holds even when $G$ is not discrete. This is the reason that the case $\mathcal{A} = PF_\Psi(G)$ is not considered in the converse of the lemma. Furthermore, it is easy to verify that the product $\square$ on $\mathcal{A}^{'} = W_\Phi(G)$ coincides with the pointwise multiplication.
\end{rem}

\begin{thm} Let $G$ be a discrete group. If $M$ is a maximal regular left/right/two-sided ideal of $\mathcal{A^{'}},$ then either $M$ contains $A_\Phi(G),$ or there exists a unique $x \in G$ such that $$M = M_x = \{\Gamma \in \mathcal{A^{'}}: \langle \Gamma,\lambda_\Psi(x)\rangle = 0 \}.$$
\end{thm}

\begin{proof}
Let $M$ be a maximal regular left/right/two-sided ideal of $\mathcal{A}^{'}$ which does not contain $A_\Phi(G).$ By Corollary \ref{s3c1}, $M \cap A_\Phi(G) = I_x = \{u \in A_\Phi(G): u(x) = 0\}$ for a unique $x \in G.$ Now choose $v \in A_\Phi(G)$ such that $v(x) \neq 0.$ As $G$ is discrete, by Lemma \ref{s3l3}, $A_\Phi(G)$ is an ideal in $\mathcal{A}^{'}$ and since $M$ is also a left/right/two-sided ideal, for any $\tilde{\Gamma} \in M,$ it follows that $v \, \square \, \tilde{\Gamma} \in M \cap A_\Phi(G).$ This implies that $(v \, \square \, \tilde{\Gamma})(x) = 0.$ Further observe that with $J(\lambda_\Psi(x))$ as the canonical image of $\lambda_\Psi(x)$ in $PM_\Psi(G)^{''},$ we have,
$$0 = (v \, \square \, \tilde{\Gamma})(x) = \langle J(\lambda_\Psi(x)),v \, \square \, \tilde{\Gamma} \rangle = \langle J(\lambda_\Psi(x)),v\rangle \, \langle J(\lambda_\Psi(x)),\tilde{\Gamma}\rangle = v(x) \, \langle \tilde{\Gamma},\lambda_\Psi(x)\rangle.$$

As $v(x) \neq 0,$ it follows that $\langle \tilde{\Gamma},\lambda_\Psi(x)\rangle = 0.$ This implies that $M \subseteq M_x,$ where $M_x = \{\Gamma \in \mathcal{A^{'}}: \langle \Gamma,\lambda_\Psi(x)\rangle = 0 \}$. Since $M_x$ is also an ideal of $\mathcal{A}^{'},$ the result follows from the maximality of $M.$
\end{proof}

Recall that a group $G$ is said to be amenable if there exists a continuous linear functional $\gamma$ on $L^\infty(G)$ such that $\| \gamma \| = \gamma(\textbf{1}) = 1$ and $\gamma(L_gf) = \gamma(f)$ for all $f \in L^\infty(G)$ and $g \in G.$ The last theorem of this section provides equivalent conditions for each of the maximal regular left ideals of $\mathcal{A}^{'}$ to contain $I_x$ for a unique $x \in G$ under the assumption that the underlying group $G$ is amenable.

\begin{thm} For an amenable group $G,$ the following are equivalent.

    $(1)$ The Banach algebra $\mathcal{A}^{'}$ has an identity contained in the left ideal $M$ generated by $A_\Phi(G).$
    
    $(2)$ Each maximal regular left ideal $M$ of $\mathcal{A}^{'}$ contains $I_x$ for a unique $x \in G.$
    
    $(3)$ The left ideal $M$ in $\mathcal{A}^{'}$ generated by $A_\Phi(G)$ contains a right identity of $\mathcal{A}^{'}.$
\end{thm}

\begin{proof}
The proof of the theorem is similar to the proof of \cite[Theorem 4.7]{GL}. It uses the fact that the given Banach algebra has a bounded approximate identity. In the case of $A_\Phi(G),$ since $G$ is amenable, \cite[Theorem 3.1]{RLSK2} gives the existence of a bounded approximate identity, and hence, the proof follows.
\end{proof}

\section{Minimal ideals} \label{sec4}

We begin this section by recalling the notion of topologically invariant mean on $PM_\Psi(G).$

Let $e$ denotes the identity of the group $G.$ A mean $\Gamma$ on $PM_\Psi(G)$ is said to be topologically invariant (\cite[Definition 6.1]{RLSK1}) if for all $T \in PM_\Psi(G)$ and $u \in A_\Phi(G),$
$$\langle \Gamma,u \cdot T\rangle = u(e) \, \langle \Gamma,T\rangle$$
holds. Moreover, it follows from \cite[Corollary 6.2]{RLSK1} that there exists a topologically invariant mean on $PM_\Psi(G).$

Motivated from \cite{DFM}, we have the following definition of topological $x$-invariance on $\mathcal{A}$ for some $x \in G.$

\begin{defn}
A linear functional $\Gamma$ on $\mathcal{A}$ is said to be topologically $x$-invariant for some $x \in G$ if
$$\langle \Gamma,u \cdot T \rangle := u(x) \, \langle \Gamma,T\rangle$$
for $u \in A_\Phi(G)$ and $T \in \mathcal{A}.$    
\end{defn}

\begin{rem}\label{s4r1}
    Since $\mathcal{A}$ is topologically invariant, we can define an $A_\Phi$-module map on $\mathcal{A}^{'}$ by 
    $$\langle u \cdot \Gamma,T\rangle := \langle \Gamma,u \cdot T\rangle$$
    for $u \in A_\Phi(G), T \in \mathcal{A}$ and $\Gamma \in \mathcal{A}^{'}.$ Also, one can verify that $u \, \square \, \Gamma = u \cdot \Gamma.$
\end{rem}

By Remark \ref{s4r1}, the equivalent definition of topological $x$-invariance for $\Gamma \in \mathcal{A}^{'}$ is
$$u \, \square \, \Gamma = u \cdot \Gamma = u(x) \, \Gamma \hspace{1cm} \text{for all} \, \, u \in A_\Phi(G).$$

For $x = e,$ this definition coincides with the definition of topological invariance on $\mathcal{A} = PM_\Psi(G).$ Further, the existence of topologically $x$-invariant functional on $\mathcal{A}$ follows from \cite[Corollary 6.2]{RLSK1}, by restriction.

For $x \in G,$ we define an action of $\lambda_\Psi(x)$ on $\mathcal{A}$ by $T \mapsto \lambda_\Psi(x)T$ for $T \in \mathcal{A}.$ The map $\lambda_\Psi(x)T: A_\Phi(G) \to \mathbb{C}$ is given by
$$\langle \lambda_\Psi(x)T,v\rangle := \langle T,L_{x^{-1}}v\rangle, \hspace{1cm} (v \in A_\Phi(G)),$$
where $L_{x^{-1}}$ denotes the left translation by $x^{-1}$ and is defined as $(L_{x^{-1}}v)(y) := v(xy),$ for $y \in G.$ 
It is easy to verify that $\lambda_\Psi(x) \mathcal{A} = \mathcal{A}.$ Using this action, for $\Gamma \in \mathcal{A}^{'}$ and $x \in G,$ one can define a map $\Gamma_x \in \mathcal{A}^{'}$ by $$\langle \Gamma_x,T\rangle := \langle \Gamma,\lambda_\Psi(x) \, T\rangle, \hspace{1cm} (T \in \mathcal{A}).$$ 

The next lemma provides a characterisation for $\Gamma \in \mathcal{A}^{'}$ to be topologically $x$-invariant in terms of $\Gamma_x.$

\begin{lem}
    Let $\Gamma \in \mathcal{A}^{'}.$ Then $\Gamma$ is topologically $x$-invariant for some $x \in G$ if and only if $\Gamma_x$ is topologically invariant.
\end{lem}

\begin{proof}
   Suppose $\Gamma$ is topologically $x$-invariant for some $x \in G.$ It is easy to verify that for $u \in A_\Phi(G)$ and $T \in \mathcal{A},$ $$\lambda_\Psi(x)(u \cdot T) = (L_x u)\cdot (\lambda_\Psi(x)T).$$
   Now,
   \begin{align*}
   \langle \Gamma_x,u \cdot T \rangle &= \langle \Gamma,\lambda_\Psi(x)(u \cdot T)\rangle\\ &= \langle \Gamma,(L_x u)\cdot (\lambda_\Psi(x)T)\rangle \\ &= (L_x u)(x) \, \langle \Gamma,\lambda_\Psi(x)T\rangle \\ &= u(e) \, \langle \Gamma_x,T\rangle,
   \end{align*}
   which implies that $\Gamma_x$ is topologically invariant.

   The converse follows on the same lines and is a consequence of the fact that
   $$u \cdot T = \lambda_\Psi(x)(L_{x^{-1}}u \cdot \lambda_\Psi(x^{-1})T)$$
   for $u \in A_\Phi(G)$ and $T \in \mathcal{A}.$
\end{proof}

\begin{defn}
    A linear functional $\Gamma$ on $\mathcal{A}$ is termed as a right annihilator of $\mathcal{A}^{'}$ if $\mathcal{A}^{'} \,\square \, \Gamma = \{0\}.$
\end{defn}

By using Remark \ref{s3r1}, \ref{s4r1} and the weak$^\ast$ continuity of the product $\square$ in $\mathcal{A}^{'}$ on the left side, one can verify that $\Gamma \in \mathcal{A}^{'}$ is a right annihilator of $\mathcal{A}^{'}$ if and only if $\Gamma$ vanishes on $A_\Phi \cdot \mathcal{A}.$

\begin{nota}
    For an element $b$ in an algebra $B,$ we denote the set $\{k \, b : k \, \, \text{is a scalar}\}$ by $\textit{\textbf{k}} \, b.$
\end{nota}

\begin{prop}\label{s4p2}
    Let $\Gamma \in \mathcal{A}^{'}.$ Assume that $\Gamma$ is either a non-zero right annihilator of $\mathcal{A}^{'},$ or topologically $x$-invariant for some $x \in G.$ Then the left ideal generated by $\Gamma$ in $\mathcal{A}^{'}$ is $\textbf{k} \, \Gamma,$ i.e., of dimension one and hence, minimal. 
\end{prop}

\begin{proof}
    The proof of the proposition follows verbatim to the proof of \cite[Proposition 5.3]{DFM}. Hence, we omit it.
\end{proof}

Here is the main theorem of this section that characterises the minimal left ideals of $\mathcal{A}^{'}.$ This is an Orlicz analogue of \cite[Theorem 5.8]{DFM}.

\begin{thm}\label{s4t6}
    Let $M$ be a left ideal in $\mathcal{A}^{'}.$ Then $M$ is minimal if and only if $M = \textbf{k} \, \Gamma,$ where $\Gamma$ is either topologically $x$-invariant for some $x \in G,$ or a non-zero right annihilator of $\mathcal{A}^{'}.$
\end{thm}

\begin{proof}
    If $M = \textit{\textbf{k}} \, \Gamma,$ where $\Gamma$ is either topologically $x$-invariant for some $x \in G$ or a non-zero right annihilator of $\mathcal{A}^{'},$ then by Proposition \ref{s4p2}, the left ideal generated by $\Gamma$ is minimal and equals $\textit{\textbf{k}} \, \Gamma.$ This implies $M$ is minimal.

    Conversely, suppose that $M$ is a left minimal ideal in $\mathcal{A}^{'}.$ Let $\Gamma \in M$ and $\Gamma \neq 0.$ If $\mathcal{A}^{'} \, \square \, \Gamma = \{0\},$ then $\Gamma$ is a right annihilator of $\mathcal{A}^{'}$ and hence, by Proposition \ref{s4p2}, $M = \textit{\textbf{k}} \, \Gamma.$ Otherwise, there exists $\Gamma_0 \in \mathcal{A}^{'}$ such that $\Gamma_0 \, \square \, \Gamma \neq 0.$ By \cite[Lemma 3.4]{filali}, $\mathcal{L}(\Gamma)$ is a maximal regular left ideal of $\mathcal{A}^{'},$ where
    $\mathcal{L}(\Gamma):= \{\tilde{\Gamma} \in \mathcal{A}^{'}: \tilde{\Gamma} \, \square \, \Gamma = 0\}.$ As the product $\square$ in $\mathcal{A}^{'}$ is weak$^\ast$ continuous on the left, it follows that $\mathcal{L}(\Gamma)$ is weak$^\ast$ closed. Since $\Gamma_0 \, \square \, \Gamma \neq 0,$ by Theorem \ref{s3t1}, there exists a unique $x \in G$ such that 
    $$\mathcal{L}(\Gamma) = \{\tilde{\Gamma} \in \mathcal{A}^{'}: \langle \tilde{\Gamma},\lambda_\Psi(x)\rangle = 0\}.$$
    Since $M$ is a minimal ideal, $M = \mathcal{A}^{'} \, \square \, \Gamma,$ i.e., $M$ is the left ideal generated by $\Gamma.$ Take $\Gamma_e \in \mathcal{A}^{'}$ such that $\Gamma_e \, \square \, \Gamma = \Gamma.$ Now, for any $v \in A_\Phi(G),$ it is easy to verify that $v - v \, \square \, \Gamma_e \in \mathcal{L}(\Gamma).$ Thus,
    $$0 = \langle v - v \, \square \, \Gamma_e,\lambda_\Psi(x)\rangle = v(x) - \langle v \cdot \Gamma_e,\lambda_\Psi(x)\rangle = v(x) - v(x) \, \langle \Gamma_e,\lambda_\Psi(x)\rangle.$$
    Choose $v \in A_\Phi(G)$ such that $v(x) = 1.$ It follows that $\Gamma_e(\lambda_\Psi(x)) = 1.$ This implies that for any $v \in A_\Phi(G),$ $v - v(x) \, \Gamma_e \in \mathcal{L}(\Gamma)$ and hence, 
    $$v \cdot \Gamma = v \, \square \, \Gamma = v(x) (\Gamma_e \, \square \, \Gamma) = v(x) \, \Gamma.$$
    This proves that $\Gamma$ is topologically $x$-invariant. The assertion is a consequence of Proposition \ref{s4p2}.
\end{proof}

The next result gives the existence of minimal idempotents in the Banach algebras $A_\Phi(G)$ and $\mathcal{B}.$ Let $\mathrm{B}$ denote any of the commutative Banach algebras $A_\Phi(G)$ or $\mathcal{B}.$ Recall that a non-zero element $b$ in a Banach algebra $B$ is said to be minimal idempotent if $b^2 = b$ and $b \, B \, b = \textit{\textbf{k}}\, b.$

\begin{rem}\label{s4r2}
    If $\varphi \in \sigma(\mathcal{B})$ and $\varphi \neq \mathcal{E}_x$ for every $x \in G,$ then $\varphi(A_\Phi(G)) = \{0\},$ where $\mathcal{E}_x(u) := u(x)$ for $u \in \mathcal{B}.$ This follows from the fact that $A_\Phi(G) \subseteq \mathcal{B}$ and the Gelfand spectrum of $A_\Phi(G)$ is homeomorphic to $G$ \cite[Corollary 3.8]{RLSK1}.
\end{rem}

\begin{thm}\label{s4p3}
    A minimal idempotent exists in $\mathrm{B}$ if and only if $G$ is discrete. 
\end{thm}

\begin{proof}
    Let $u$ be a minimal idempotent in $\mathrm{B}.$ Then $u^2 = u$ and $u \, \mathrm{B} \, u = \textit{\textbf{k}}\, u.$ As $\mathrm{B}$ is commutative, for any $v \in \mathrm{B},$ $v \, u = \varphi(v) \, u,$ where $\varphi(v)$ is a scalar. Consider the linear map from $\mathrm{B}$ to $\mathbb{C},$ given by $v \mapsto \varphi(v).$ Observe that for $v_1, v_2 \in \mathrm{B},$
    $$\varphi(v_1 \, v_2) \, u = v_1 \, v_2 \, u = (v_1 \, u) \, (v_2 \, u) = \varphi(v_1) \, \varphi(v_2) \, u.$$
    Since $u \neq 0,$ it follows that $\varphi(v_1 \, v_2) = \varphi(v_1) \, \varphi(v_2),$ i.e., $\varphi \in \sigma(\mathrm{B}).$ Again, as $u \neq 0$ and $A_\Phi(G) \subseteq \mathrm{B},$ we have,
    $\varphi(A_\Phi(G)) \neq \{0\}.$ By Remark \ref{s4r2}, $\varphi = \mathcal{E}_x$ for some $x \in G.$ Observe that $v \, u = \varphi(v) \, u = \mathcal{E}_x(v) \, u = v(x) \, u$ for any $v \in \mathrm{B}.$ Our claim is that $u = \chi_{\{x\}}.$ Let $y \in G$ such that $y \neq x.$ Choose $v \in A_\Phi(G) \subseteq \mathrm{B}$ such that $v(x) \neq v(y).$ Now, $v(y) \, u(y) = (v \, u)(y) = (v(x) \, u)(y) = v(x) \, u(y)$ which implies that $u(y) = 0.$ Thus, $u = \chi_{\{x\}} \in \mathrm{B} \subseteq \mathcal{C}(G)$ and hence, $G$ is discrete. The converse is obvious.
\end{proof}

Recall that an algebra $B$ is said to be semi-prime if $\{0\}$ is the only bi-ideal $M$ of $B$ with $M^2 = \{0\}.$ Since $\mathrm{B}$ is semi-simple, by \cite[Proposition 5, Pg. 155]{BD}, it follows that $\mathrm{B}$ is semi-prime. The last corollary is a direct consequence of Theorem \ref{s4p3} and \cite[Proposition 6 (i), Pg. 155]{BD}.

\begin{cor}
    The minimal ideals $M$ exist in $\mathrm{B}$ if and only if $G$ is discrete and $M = \textbf{k} \, \chi_{\{x\}}$ for some $x \in G.$
\end{cor}

\section*{Acknowledgement}
The first author is grateful to the Indian Institute of Technology Delhi for the Institute Assistantship.

\section*{Data Availability} 
Data sharing does not apply to this article as no datasets were generated or analysed during the current study.

\section*{Competing Interests}
The authors declare that they have no competing interests.

\bibliographystyle{acm}
\bibliography{article3}

\end{document}